\documentclass[10pt]{amsart}

\usepackage{amssymb, verbatim}

\theoremstyle{plain}
\newtheorem{theorem}{Theorem}[section]
\newtheorem{proposition}[theorem]{Proposition}
\newtheorem{corollary}[theorem]{Corollary}

\theoremstyle{definition}
\newtheorem{definition}[theorem]{Definition}

\newtheorem{remark}[theorem]{Remark}

\def\accute{\'}
\def\<{\langle}
\def\>{\rangle}
\def\geqs{\geqslant}
\def\leqs{\leqslant}
\def\a{\alpha}
\def\b{\beta}
\def\d{\delta}
\def\k{\kappa}
\def\v{\varphi}
\def\w{\omega}
\def\H{\mathcal H}
\def\V{\mathcal V}
\def\/{\kern 0.05em}

\DeclareMathOperator{\ricci}{Ricci} 

\begin{document}
\title{The geometry of generalised Cheeger-Gromoll metrics}

\author{M. Benyounes}
\address{D{\accute e}partement de Math{\accute e}matiques \\
Universit{\accute e} de Bretagne Occidentale \\
6, avenue Victor Le Gorgeu \\
CS 93837, 29238 Brest Cedex 3\\
France} \email{Michele.Benyounes@univ-brest.fr}

\author{E. Loubeau}
\address{D{\accute e}partement de Math{\accute e}matiques \\
Universit{\accute e} de Bretagne Occidentale \\
6, avenue Victor Le Gorgeu \\
CS 93837, 29238 Brest Cedex 3\\
France} 
\email{Eric.Loubeau@univ-brest.fr}

\author{C.~M. Wood}
\address{Department of Mathematics \\
University of York \\
Heslington, York Y010 5DD\\
U.K.} \email{cmw4@york.ac.uk}



\date{\today}

\begin{abstract}
We study the geometry of the tangent bundle equipped with a two-parameter family of metrics, deforming the Sasaki and Cheeger-Gromoll metrics. After deriving the expression for the Levi-Civita connection, we compute the Riemann curvature tensor and the sectional, Ricci and scalar curvatures.  We identify all metrics whose restriction to the fibres has positive sectional curvature.
When the base manifold is a space form, we characterise metrics with non-negative sectional curvature and show that one can always find parameters ensuring positive scalar curvature.  This extends to compact manifolds and, under some curvature conditions, to general manifolds.
\end{abstract}

\maketitle

\section{Introduction}

Let $(M,g)$ be a smooth Riemannian manifold of dimension $n$, with tangent bundle $\pi\colon TM\to M$, and let $K\colon TTM\to TM$ denote the connection map for the Levi-Civita connection~\cite{Dom62}.  Then the horizontal distribution $\H=\ker(K)$ is complementary to the vertical distribution $\V=\ker(d\pi)$:
$$
TTM=\H\oplus\V,
$$
and the metric $g$ of $M$ may be lifted to a metric $h$ on $TM$:
$$ 
h(A,B) = \< d\pi(A),d\pi(B) \> + \<KA,KB\>,
$$
for all $A,B \in T_eTM $ and all $e\in TM$, where we have abbreviated $g(X,Y)=\<X,Y\>$.  This is the well-known {\sl Sasaki metric\/} \cite{Sas58}.  It has been extensively studied and shown to be rather rigid, especially when it comes to the existence of harmonic sections of $TM$, all of which are parallel when $M$ is compact \cite{Ish79, M-T, Nou77}, and hence trivial if the Euler characteristic of $M$ is non-zero.  The rigidity of $h$ may be overcome to a certain extent by restricting it to the sphere subbundle $SM(r)$ of radius $r>0$, and considering the harmonic section variational problem constrained to sections of length $r$.  This has been very fruitful (cf.~the survey article \cite{Gil05}), but is limited to manifolds $M$ whose Euler characteristic vanishes.  A more recent approach has been to introduce new metrics on $TM$, which satisfy the natural conditions that $\pi$ is a Riemannian submersion (so that the geometry of $M$ is reflected in that of $TM$) and $\H$ and $\V$ are orthogonal, with $\V$ totally geodesic, but whose restriction to the tangent spaces is non-Euclidean, and which in some cases may be regarded as vertical ``geometric compactifications'' of $TM$.   This has led to the $2$-parameter family of {\sl generalised Cheeger-Gromoll\/} metrics~\cite{BLW1}:
$$
h_{p,q}(A,B) 
= \< d\pi(A),d\pi(B) \> + \w^p \big[ \< KA,KB \> + q \<KA,e\> \<KB,e\> \big], 
$$
for all $p,q\in\mathbb R$, where $\w(e)=(1+|e|^2)^{-1}$, with $|e|^2=\<e,e\>$.  The original Cheeger-Gromoll metric corresponds to $p=q=1$ (\cite{C-G72, M-T}), and the Sasaki metric to $p=q=0$.  If $q\geqs0$ then $h_{p,q}$ is a {\it bona fide\/} Riemannian metric on $TM$, for all $p$.  However if $q<0$ then $h_{p,q}$ is positive definite on the (open) ball subbundle $BM(r)$ of radius $r=1/\sqrt{-q}$, which we refer to as the {\sl Riemannian ball bundle\/} of $h_{p,q}$.  The restriction of $h_{p,q}$ to the sphere bundle $SM(r)$ bounding the Riemannian ball bundle is also positive definite, for all $p$, and this is true for all the sphere subbundles of $TM$;  however the {\sl canonical vertical vector\/} in $TTM$ (see below) along $SM(r)$ is $h_{p,q}$-null, and $h_{p,q}$ is Lorentzian on the interior of the complement of the Riemannian ball bundle.  This family of metrics is proving to be very interesting in the theory of harmonic maps/sections, exhibiting a delicate balance between flexibility and rigidity (cf.~\cite{BLW1,BLW2}).  For example, when suitably scaled, the vector field tangent to standard Morse flow on the round sphere $S^n$ ($n\geqs3$) is a harmonic section of $TS^n$ with respect to a unique metric $h_{p,q}$.  On the other hand, when $M$ is compact of non-zero Euler characteristic there is a substantial ``Liouville region'' of the $(p,q)$-plane, where the only harmonic section with respect to $h_{p,q}$ is the zero section, and when $M$ is non-compact there is an analogous ``Bernstein region''.  The Cheeger-Gromoll and Sasaki metrics ``belong'' to both of these regions. However, in this paper we concentrate on the geometry of $h_{p,q}$, which is of interest in its own right. 

In \S2 we determine via standard methods the Levi-Civita connection of $h_{p,q}$, and compute the Riemann curvature tensor and sectional, Ricci and scalar curvatures.  We note that when $q<0$ all these curvatures are unbounded on the Riemannian ball bundle, but if $p+q=1$ the sectional curvatures of $2$-planes tangent to the boundary remain bounded (Corollary \ref{cor1}), making these the ``most regular'' metrics with $q<0$.
We also deduce that, in order for $h_{p,q}$ to be an Einstein metric, $(M,g)$ must be Einstein with harmonic curvature; in this respect, the $h_{p,q}$ behave similarly to the Sasaki metric.
However, our main result is Theorem \ref{prop3a}, which gives some precision to the sense in which $h_{p,q}$ geometrically compactifies the tangent spaces, by identifying the region of the $(p,q)$-plane such that the metric induced by $h_{p,q}$ on the fibres of $TM$ (or, when $q<0$, of the Riemannian ball bundle) has positive sectional curvature.  This region lies entirely in the half plane $2p+q>0$; when $n\geqs3$ it is a region $\Gamma$ independent of the dimension of $M$, but when $n=2$ it is a considerably larger region $\Gamma'$.  In all dimensions, it is not an open subset, but it is connected, and for each $q\in\Bbb R$ contains at least one point $(p,q)$, and conversely for each $p>-8$.  Since it contains the parameters $(p,q)=(1,1)$ of the Cheeger-Gromoll metric, it is not complementary to the Liouville or Bernstein regions of the $(p,q)$-plane.  Another notable element of both $\Gamma$ and $\Gamma'$ is $(p,q)=(2,0)$, which are the coordinates of the stereographic compactification of the fibres of $TM$.  Interestingly, when $n\geqs3$ the only $(p,q)\in\Gamma$ with $q<0$ are those with $p+q=1$.

In \S3 we study the curvature of $h_{p,q}$ when the base manifold has constant sectional curvature $c$.  After simplifying the curvature expressions (Proposition \ref{prop6}), we determine the region of the $(p,q)$-plane for which $h_{p,q}$ has non-negative sectional curvature (Theorem \ref{prop7}).  This region is non-empty if and only if $c\geqs0$.  If $n\geqs3$ then it is a region $\Delta_c$ independent of the dimension of $M$; if $n=2$ then it is a considerably larger region $\Delta_c'$.  In all dimensions, if $c\leqs 4/3$ then it is independent of $c$, whereas if $c>16/3$ then it lies entirely in the lower half of the $(p,q)$-plane, with $q<0$ if $n\geqs3$, and $q\leqs0$ if $n=2$.  
The sectional curvature of $h_{p,q}$ can never be (strictly) positive, and in the remainder of the paper we examine the weaker condition of positive scalar curvature.  A necessary condition for the scalar curvature $\tilde s=\tilde s_{p,q}$ of $h_{p,q}$ to be positive is $2p+q>-c$, which no longer precludes the possibility of $c<0$.  If $(M,g)$ is flat then $\tilde s>0$ precisely when the fibres of $TM$ have positive sectional curvature, as determined in Theorem \ref{prop3a}.  In Theorem \ref{prop8} we determine subregions of $\Gamma$ and $\Gamma'$ for which $h_{p,q}$ has positive scalar curvature.  This paves the way for Theorem \ref{thm1}, which states that for any space form there exist parameters $(p,q)$ such that $\tilde s_{p,q}>0$.  When $c\neq0$ and $n\geqs3$, all the metrics $h_{p,q}$ of Theorem \ref{thm1} have $q<0$, and consequently are only defined on the Riemannian ball subbundle, with no possibility of extension to $TM$.  This problem is rectified by Theorem \ref{thm3}, which shows that for every space form $M$ it is possible to find generalised Cheeger-Gromoll metrics on $TM$ of positive scalar curvature.  When $n\geqs3$ and $c\neq0$, the sectional curvatures of the tangent spaces with respect to these metrics are in general no longer entirely positive (or non-negative).  Both Theorems \ref{thm1} and \ref{thm3} extend to compact manifolds (Theorem \ref{thm2}), and general Riemannian manifolds under some curvature conditions (Theorem \ref{prop9}).  These results generalise work of Gudmundson and Kappos~\cite{Gud-Kap1}, and Sekizawa~\cite{Sek91}, and are particularly interesting in view of \cite[Proposition 5.4]{BLW2}, yielding new examples of harmonic maps from compact manifolds into non-flat manifolds of non-negative sectional curvature.

\newpage
\section{Generalised Cheeger-Gromoll geometry.}

Let $M$ be a Riemannian $n$-manifold, not necessarily compact, or even complete.  The generalised Cheeger-Gromoll metric $h_{p,q}$ on $TM$ 
can alternatively be described in terms of horizontal and vertical lifts of tangent vectors to $M$ (cf.~\cite{Dom62}):
\begin{align*}
h_{p,q}(X^h,Y^h) &=\<X,Y\>, \\
h_{p,q}(X^h,Y^v) &=0, \\
h_{p,q}(X^v,Y^v) &=\w^p(\<X,Y\>+q\<X,e\>\<Y,e\>),
\end{align*}
where $X,Y \in T_{\pi(e)}M$.  Most of our results will be expressed in this way.

\begin{remark}\label{rem1}
Various classes of metrics have been introduced on the tangent bundle. For example, generalised Cheeger-Gromoll metrics lie within the set of {\sl Kaluza-Klein metrics\/} studied in~\cite{Wood,Wood1}, for which the bundle projection $\pi$ is a Riemannian submersion with totally geodesic fibres, and for such metrics an Eells-Sampson type existence theorem for harmonic sections of fibre bundles is proved (although this requires the bundle to have compact, negatively curved fibres).  These Kaluza-Klein metrics are, in turn, special cases of {\sl $g$-natural metrics\/} introduced in~\cite{KS} (cf. also~\cite{ab1,ab2}) which are characterised, at a point $e \in TM$, by the relations:
\begin{align*}
G(X^h , Y^h)_e
&= A(|e|^2) g_{\pi(e)}(X,Y) + B(|e|^2)g_{\pi(e)}(X,e)g_{\pi(e)}(Y,e), \\
G(X^h , Y^v)_e
&= C(|e|^2) g_{\pi(e)}(X,Y) + D(|e|^2)g_{\pi(e)}(X,e)g_{\pi(e)}(Y,e), \\
G(X^v , Y^v)_e
&= E(|e|^2) g_{\pi(e)}(X,Y) + F(|e|^2)g_{\pi(e)}(X,e)g_{\pi(e)}(Y,e),
\end{align*}
where $A,B,C,D,E$ and $F$ are positive functions.  Note that, in general, for $g$-natural metrics,  the bundle projection is no longer a Riemannian submersion, and the fibres not totally geodesic.
\end{remark}

By standard properties of  metrics on $TM$ such that $\H$ and $\V$ are orthogonal and $\pi$ is Riemannian submersion~\cite{ON66}, and using the Koszul formula, we compute the Levi-Civita connection of $h_{p,q}$.  It is convenient to define:
$$
BM_q=\{e\in TM:q|e|^2>-1\},
\qquad
SM_q=\{e\in TM:q|e|^2=-1\}.
$$
Notice that if $q\geqs0$ then $BM_q=TM$ and $SM_q=\emptyset$, and if $q<0$ then $BM_q$ is the Riemannian ball bundle of $h_{p,q}$ and $SM_q$ is its sphere bundle boundary.  We now introduce the following smooth $1$-parameter deformation $\w_q$ of $\w$:
$$
\w_q(e)=\frac{1}{1+q|e|^2},
$$
defined for all $q\in\Bbb R$ and all $e\in BM_q$.  Finally let $U$ denote the {\sl canonical vertical vector field\/} on $TM$: $U(e)=e^v$ for all $e\in TM$.  Note that $U$ is not a unit vector field with respect to any metric $h_{p,q}$; indeed, if $q<0$ then the restriction of $U$ to $SM_q$ is $h_{p,q}$-null for all $p$.  The results of our computations are now summarised as follows.

\begin{proposition}\label{prop1}
Let $(M,g)$ be a Riemannian manifold, and let $h_{p,q}$ be a generalised Cheeger-Gromoll metric on $BM_q$.  Let $\nabla$ (respectively, $R$) denote the Levi-Civita connection (respectively, Riemann tensor) of $M$.  Then the Levi-Civita connection $\widetilde\nabla$ of $h_{p,q}$ satisfies, at $e\in BM_q$:
\begin{align}
\widetilde\nabla_{X^h}Y^h
&= (\nabla_X Y)^h -\tfrac12 [R(X,Y)e]^v, \\
\widetilde\nabla_{X^h}Y^v
&= (\nabla_X Y)^v +\tfrac12 \w^p [R(e,Y)X]^h, \\
\widetilde\nabla_{X^v}Y^h
&= \tfrac12 \w^p [R(e,X)Y]^h, \\
\widetilde\nabla_{X^v}Y^v
&= \w_q \big[ (p\w+q) \<X,Y\>  + pq\w \<X,e\> \<Y,e\> \big] U 
- p\/\w \big[\<X,e\>Y + \<Y,e\>X\big]^v, \label{eq4}
\end{align}
for all $X\in T_{\pi(e)}M$ and $Y\in C^{\infty}(TM)$ {\rm(}the set of smooth vector fields on $M${\rm)}.
\end{proposition}

\noindent
{\bf Note.}
Our convention for the Riemann curvature tensor is:
$$
R(X,Y)=[\nabla_X,\nabla_Y]-\nabla_{[X,Y]}.
$$
\vskip1ex
Comparison with \cite{Kow71} shows that equations (1)--(3) are remarkably similar to the Sasaki case, to which they reduce when $p=0$, for all $q$.  The major difference is equation (4), from which it follows that $\widetilde\nabla$ is in general singular on $SM_q$ when $q<0$.  (Note that the restriction of $\widetilde\nabla$ to $SM_q$ is not the Levi-Civita connection of the restriction of $h_{p,q}$; the latter is in fact a smooth Riemannian metric with respect to which the fibres of $SM_q$ are spheres of radius $\w^{-p}$.)

Lengthy but straightforward computations involving Proposition~\ref{prop1}, the Dombrowski Lie bracket formula~\cite{Dom62}, and Gudmundsson-Kappos' expressions for covariant differentiation of the vertical and horizontal lifts of a bundle endomorphism~\cite{Gud-Kap2}, determine the Riemann curvature tensor of $h_{p,q}$.  It is convenient to introduce the following curvature-type tensor $r$ on $M$:
\begin{equation}
r(X,Y)Z=\<Y,Z\>X-\<X,Z\>Y.
\label{littler}
\end{equation}
Again we summarise the results.

\newpage
\begin{proposition}\label{prop2}
The curvature tensor $\tilde R$ of $h_{p,q}$ is given, at $e\in BM_q$, by:
\begin{align*}
&({\rm i})\quad  
\tilde R(X^h,Y^h)Z^h 
= [R(X,Y)Z]^h + \tfrac12 \big[ (\nabla_Z R)(X,Y)e \big]^v \\
&\qquad\qquad
- \tfrac14\w^p\big[ R(e,R(Y,Z)e)X - R(e,R(X,Z)e)Y 
- 2 R(e,R(X,Y)e)Z\big]^h, \\
&({\rm ii})\quad 
\tilde R(X^h,Y^h)Z^v 
= [R(X,Y)Z]^v + \tfrac12 \w^{p} \big[
(\nabla_X R)(e,Z)Y - (\nabla_Y R)(e,Z)X\big]^h \\
&\qquad\qquad
+\tfrac14 \w^p \big[ R(Y,R(e,Z)X)e - R(X,R(e,Z)Y)e \big]^v
- p\/\w \<Z,e\> [R(X,Y)e]^v \\
&\qquad\qquad\qquad 
+ (p\/\w + q)\w_q \< R(X,Y)e,Z\> U, \\
&({\rm iii})\quad 
\tilde R(X^h,Y^v)Z^h 
= \tfrac12\w^p \big[(\nabla_{X}R)(e,Y)Z\big]^h 
- \tfrac14\w^p \big[R(X,R(e,Y)Z)e\big]^v \\
&\qquad\qquad
-\tfrac{p}{2}\/\w \<Y,e\> [R(X,Z)e]^v + \tfrac12 [R(X,Z)Y]^v 
+ \tfrac12(p\/\w +q)\w_q \< R(X,Z)e,Y \> U, \\
&({\rm iv})\quad 
\tilde R(X^h,Y^v)Z^v 
= \tfrac{p}{2}\/\w^{p+1} \big[ \<Y,e\> R(e,Z)X - \<Z,e\> R(e,Y)X\big]^h \\
&\qquad\qquad
- \tfrac12 \w^p [R(Y,Z)X]^h 
- \tfrac14 \w^{2p} \big[R(e,Y)R(e,Z)X\big]^h, \\
&({\rm v})\quad
\tilde R(X^v,Y^v)Z^h 
= \w^p [R(X,Y)Z]^h 
+ p\/\w^{p+1} \big[ \<Y,e\> R(e,X)Z - \<X,e\> R(e,Y)Z \big]^h \\
&\qquad\qquad
+ \tfrac14 \w^{2p} [R(e,X)R(e,Y)Z - R(e,Y)R(e,X)Z]^h, \\
&({\rm vi})\quad 
\tilde R(X^v,Y^v)Z^v 
= A\,\<Z,e\> [ r(X,Y)e ]^v 
+ B [ r(X,Y)Z ]^v 
+ C \<r(X,Y)Z,e\> U,
\end{align*}
where:
\begin{align*}
A &= p\/\w\w_q ((p+2q-2)\w -q), \\
B &= \w_q(p^2\w^2 - p(p-2)\w + q), \\
C &= \w_q^2 \big( p(p-2)(1-q)\w^2 + pq(p-3)\w - q^2 \big).
\end{align*}
\end{proposition}

\noindent
{\bf Note.}
$A$, $B$ and $C$ are related by:
\begin{equation}
\label{eqABC}
\w_q(A-qB)=C.
\end{equation}

The expressions for the Levi-Civita connection and Riemann curvature tensor of $h_{p,q}$ given in Propositions \ref{prop1} and \ref{prop2} could also be recovered as special cases of the formulas found in~\cite{ab1} and~\cite{ab2}.

The most significant manifestations of the Riemann tensor are of course the sectional, Ricci and scalar curvatures.

\begin{proposition}\label{prop3}
Let $K$ denote the sectional curvature of $(M,g)$.  Then the sectional curvature $\tilde K$ of $h_{p,q}$ satisfies the following, at $e\in BM_q$:
\begin{align}
\tilde K(X^h \wedge Y^h)
&=K(X\wedge Y)-\tfrac34 \w^p |R(X,Y)e|^2, \label{eq6}\\
\tilde K(X^h \wedge Y^v)
&=\frac{\w^p}{4(1+q\<Y,e\>^2 )} |R(e,Y)X|^2, \label{eq7}\\
\tilde K(X^v\wedge Y^v)
&=\frac{\w^{-p}}{1+q(\<X,e\>^2 + \<Y,e\>^2)}
\big[ A (\<X,e\>^2 + \<Y,e\>^2) + B\big], \label{eq8}
\end{align}
for all orthonormal vectors $X$ and $Y$ of $T_{\pi(e)} M$.  The functions $A$ and $B$ on $BM_q$ are as defined in Proposition~\ref{prop2}.
\end{proposition}

\begin{proof}
Proposition~\ref{prop2}~(i) yields:
\begin{align*}
h_{p,q}(\tilde R(X^h,Y^h)Y^h,X^h) 
&= \< R(X,Y)Y, X\> + \tfrac34 \w^p \< R(e, R(X,Y)e)Y, X\> \\
&= K(X \wedge Y)-\tfrac34 \w^p | R(X, Y)e|^2 ,
\end{align*}
and formula~\eqref{eq6} then follows because:
$$
|X^h \wedge Y^h|_{p,q}^2 = |X|^2 |Y|^2 - \<X,Y\>^2 = 1 .
$$
To obtain \eqref{eq7}, first observe that
$$
|X^h \wedge Y^v |_{p,q}^2 = \w^p(1 + q \< Y, e \>^2 ).
$$
Then by  Proposition~\ref{prop2}~(iv) we have:
$$
h_{p,q}(\tilde R(X^h,Y^v)Y^v,X^h) 
= -\tfrac14 \w^{2p} \< R(e,Y)R(e,Y)X,X\>  
= \tfrac14 \w^{2p} |R(e,Y)X|^2,
$$
and we  deduce \eqref{eq7}. Finally, since 
$$
|X^v \wedge Y^v|_{p,q}^2 
= \w^{2p} \big[ 1+q\<X,e\>^2 +q\<Y,e\>^2 \big]
$$
and
$$
\tilde R(X^v,Y^v)Y^v 
= (A\,\<Y,e\>^2 +B) X^v - A\,\<X,e\>\<Y,e\>Y^v +C\<X,e\> U,
$$
we obtain:
\begin{align*}
h_{p,q}(\tilde R(X^v,Y^v)Y^v,X^v)
&= ( A\,\<Y,e\>^2 +B) h_{p,q}(X^v,X^v) \\
&- A\,\<X,e\>\<Y,e\> h_{p,q}(X^v,Y^v) + C\<X,e\> h_{p,q}(X^v,U) \\
&= \w^p \big[ A(\<Y,e\>^2 + \<X,e\>^2) + B \big]
\end{align*}
from \eqref{eqABC}, which yields the expression for $\tilde K(X^v\wedge Y^v)$.
\end{proof}

We will say that a generalised Cheeger-Gromoll metric has {\sl flat fibres\/} if its restriction to each fibre of the Riemannian ball bundle of $M$ is flat.  The Sasaki metric has flat fibres, and is flat if and only if the base manifold is flat.  In fact, these are the only flat generalised Cheeger-Gromoll metrics.

\begin{corollary}\label{cor2}
The only generalised Cheeger-Gromoll metric with flat fibres is the Sasaki metric.
\end{corollary}

\begin{proof}
Let $e\in BM_q$, and let $\Pi$ be a vertical $2$-plane in $T_eTM$ which contains $U(e)$.  (If $n=2$ then $\Pi$ is unique, and is the only vertical $2$-plane.)  Since the tangent spaces of $M$ are totally geodesic, it suffices to consider $\tilde K(\Pi)$, and it follows from \eqref{eq8} that:
\begin{align*}
\tilde K(\Pi)
&= \w^{-p}\w_q(A|e|^2+B) \\
&= \w^{-p}\w_q^2\big( 2p(1-q)\w^2 + 3pq\w + q(1-p) \big) \\
&= \w^{2-p}\w_q^2 P(|e|^2),
\end{align*}
where:
\begin{equation}
\label{eqP}
P(t)=P_{p,q}(t)=2p+q+(p+2)q\/t+(1-p)q\/t^2.
\end{equation}
The only parameters for which $P(t)$ vanishes identically are $p=q=0$.
\end{proof}

If $q<0$ then in general $A$ and $B$ are singular on $SM_q$, and therefore by equation \eqref{eq8} the sectional curvature is undefined for $2$-planes in $T_eTM$ for all $e\in SM_q$.  There is however an exceptional case.
It is convenient to introduce the following function:
\begin{equation}
\label{eqmu}
\mu(p)=\frac{p^p}{(p-1)^{p-1}},
\end{equation}
which is defined for $p\geqs 1$, strictly increasing with $\mu(1)=1$, and unbounded above.

\begin{corollary}\label{cor1}
If $q<0$ then $\tilde K(\Pi)$ is non-singular for all $2$-planes $\Pi$ tangent to $SM_q$ if and only if $p+q=1$.  In particular, if $M$ has dimension $n\geqs3$ and $\Pi$ is vertical then $\tilde K(\Pi)=\mu(p)$.
\end{corollary}

\begin{proof}
For all $e\in SM_q$ we have:
\begin{align*}
T_e(SM_q)
&=\{A\in T_eTM:\<KA,e\>=0\}  \\
&=\H_e\oplus\{Y^v:Y\in T_{\pi(e)}M\;\text{and}\; \<Y,e\>=0\}.
\end{align*}
It follows from equation \eqref{eq7} that if $\Pi=X^h\wedge Y^v$ then:
$$
\tilde K(X^h\wedge Y^v)=\tfrac14\/\w^p\,|R(e,Y)X|^2.
$$
Now if $p+q=1$ then $A$ (respectively, $B$) extends to the smooth function $(q-1)\w^2$ (respectively, $(1-q)\w^2+w$) on $TM$.  In particular, the value of $B$ on $SM_q$ is $-q$, and it therefore follows from equation \eqref{eq8} that if $\Pi$ is vertical then: 
$$
\tilde K(\Pi)
=-q\w^{-p}
=-q(1-1/q)^p
=(p-1)\left(\frac{p}{p-1}\right)^p=\mu(p).
$$
If $p+q\neq1$ then $A$ and $B$ are unbounded on $BM_q$.
\end{proof}

If $q<0$ and $p+q=1$, then $\tilde K$ is in fact unbounded on $BM_q$.  For, if $\Pi$ is a vertical $2$-plane in $T(BM_q)$ containing the canonical vertical vector $U$, then, noting that $P_{p,q}(t)$ has a root at $t=-1/q$ when $p+q=1$, it follows from the proof of Corollary \ref{cor2} that:
$$
(1+t)^{1+q}\tilde K(\Pi)
=\frac{2-q+qt}{1+qt},
$$
where $t=|e|^2$.  It follows that the sectional curvature of $h_{p,q}$ is singular on the boundary of the Riemannian ball bundle, for all $(p,q)$ with $q<0$.  Note also, in passing, that $\tilde K(\Pi)>0$, the significance of which will be seen in Theorem \ref{prop3a}.

Setting $e=0$ in Proposition \ref{prop3} yields:
$$
\tilde K(X^h\wedge Y^h)=K(X\wedge Y),
\qquad
\tilde K(X^h\wedge Y^v)=0,
\qquad
\tilde K(X^v\wedge Y^v)=B=2p+q.
$$
It follows that the sectional curvature of $h_{p,q}$ is never strictly positive (or negative), and necessary conditions for $h_{p,q}$ to have non-negative sectional curvature are $K\geqs0$ and $2p+q\geqs0$.  In our next result we refine the condition $2p+q>0$ which is necessary for the fibres of $BM_q$ to have strictly positive sectional curvature.
It is convenient to introduce the following function:
\begin{equation}
\label{eqlambda}
\lambda(p)=
\dfrac{8(1-p)}{8+p},
\quad
p\neq-8,
\end{equation}
which parametrises the hyperbola in the $(p,q)$-plane with equation $pq+8p+8q=8$.  We now define subsets $\Gamma$ and $\Gamma'$ of the half plane $2p+q>0$ as follows.

\begin{definition}\label{def1}
Let $\Gamma_+^i$ ($i=1,2,3$) be the following region of the $(p,q)$-plane:
\begin{align*}
\Gamma_+^1
&= \{(p,q): -8<p\leqs -2 \mbox{ and } q>\lambda(p)\}, \\
\Gamma_+^2
&= \{(p,q): -2\leqs p\leqs 0 \mbox { and } 2p+q>0\}, \\ 
\Gamma_+^3 
&= \{(p,q): 0\leqs p\leqs 1 \text{ and } q>0\},
\end{align*}
and set:
$$
\Gamma_+=\Gamma^1_+\cup\Gamma^2_+\cup\Gamma^3_+.
$$
Define further regions:
\begin{align*}
\Gamma_-&=\{(p,q):p+q=1 \text{ and } q<0\}, \\
\Gamma_-'&=\{(p,q):p+q\geqs1 \text{ and } q<0\}, \\
\Gamma_Z&=\{(p,0):0<p\leqs2\}, \\
\Gamma_Z'&=\{(p,0):p>0\},
\end{align*}
and set: 
$$
\Gamma=\Gamma_-\cup\Gamma_Z\cup\Gamma_+,
\qquad 
\Gamma'=\Gamma_-'\cup\Gamma_Z'\cup\Gamma_+.
$$  
Note that $\Gamma\subset\Gamma'$, and $\Gamma'$ is a subset of the half plane $2p+q>0$.  The point $(2,0)\in\Gamma$ parametrises the metric whose vertical component is (upto homothety) the stereographic metric on $\Bbb R^n$.
\end{definition}

It is convenient to define, for each $q\in\mathbb{R}$, the following interval:
\begin{equation}
\label{eqIq}
I_q=\{t>0:qt>-1\}.
\end{equation}
If $q\geqs0$ then $I_q=\mathbb{R}^+$, whereas if $q<0$ then $I_q=[0,-1/q)$.
  
\begin{theorem}\label{prop3a}
Let $(M,g)$ be any Riemannian $n$-manifold.  If $n\geqs3$ {\rm(}respectively, $n=2${\rm)} then $\tilde K(\Pi)>0$ for all vertical $2$-planes $\Pi$ in $T(BM_q)$ precisely when  $(p,q)\in\Gamma$ {\rm(}respectively, $(p,q)\in\Gamma'${\rm)}.
In particular, the regions $\Gamma$ and $\Gamma'$ characterise the values of $(p,q)$ for which the metric induced by $h_{p,q}$ on the fibres of $BM_q$ has positive sectional curvature.

\end{theorem}

\begin{proof}
We assume throughout that $2p+q>0$.

Suppose first that $n=2$.  From the proof of Corollary \ref{cor2}, the sign of $\tilde K(\Pi)$ is controlled by the sign over $I_q$ of the polynomial $P(t)$ defined in \eqref{eqP}.  
If $q=0$ then $P(t)=2p$, so $P(I_q)>0$ for all $p>0$.
If $p=1$ then $P(t)=3qt+q+2$, so $P(I_q)>0$ for all $q\geqs0$, but $P(-1/q)=q-1<0$ for all $q<0$.  
In general, when $p\neq1$ and $q\neq0$, $P(t)$ is quadratic, with discriminant: 
$$
D=pq(pq+8p+8q-8),
$$ 
and critical point $t_0=(p+2)/2(p-1)$, which satisfies $t_0<0$ if and only if $-2<p<1$.  When $D>0$, let $t_+$ (respectively, $t_-$) denote the maximum (respectively, minimum) root of $P(t)$.
Suppose first that $q>0$.  If $p>1$ then $P(t)$ is eventually negative.  If $0\leqs p<1$ then $P(\Bbb R^+)>0$.  If $p<0$ then $P(t_0)$ is a global minimum.  When $-2<p<0$ we have $t_0<0$, so $P(\Bbb R^+)>0$ precisely when $2p+q>0$; when $p\leqs -2$ we have $t_0\geqs0$, so $P(\Bbb R^+)>0$ precisely when $D<0$, which is the case if and only if $q>\lambda(p)$ and $p>-8$.  In summary, $P(I_q)>0$ for $q>0$ if and only if $(p,q)\in\Gamma_+$.
Suppose now that $q<0$.  It suffices to consider $p\geqs0$.  If $0\leqs p<1$ then $P(t_0)$ is a global maximum, with $t_0<0$ and $D>0$.  Hence $P(I_q)>0$ precisely when $t_+\geqs-1/q$:
$$
\sqrt D\geqs 2-2p-2q-pq,
$$
which rearranges to:
$$
(p-1)(q-1)(p+q-1)\geqs 0,
$$
and is therefore impossible. 
If $p>1$ then $P(t_0)$ is a global minimum, with $t_0>0$.  If $q>\lambda(p)$ then $D<0$, so $P(I_q)>0$ for all $q<0$.  If $q=\lambda(p)$ then $D=0$, so $P(I_q)>0$ precisely when $t_0\geqs-1/q$, which after some rearrangement is equivalent to $p+q\geqs1$.  If $q<\lambda(p)$ then $D>0$, so $P(I_q)>0$ precisely when $t_-\geqs-1/q$:
$$
\sqrt D\leqs 2p-2-2q-pq,
$$
which rearranges to:
$$
(p-1)(q-1)(p+q-1) \leqs 0,
$$
and is therefore again equivalent to $p+q\geqs1$.
In summary, $P(I_q)>0$ for $q<0$ if and only if $(p,q)\in\Gamma_-$.

Now suppose that $n\geqs3$.
For any $e\in BM_q$, choosing $X,Y$ in \eqref{eq8} such that $e$ lies in the plane spanned by $X$ and $Y$, it follows that the conditions for $\tilde K(\Pi)\geqs0$ when $n=2$ are all necessary when $n\geqs3$.  Additional necessary conditions may be obtained by inspecting the sign of $\tilde K(\Pi)$ when $\Pi$ projects to a $2$-plane in $T_{\pi(e)}M$ orthogonal to $e$, and by \eqref{eq8} this is determined by $B$:
$$
B=\w_q\big(p\w^2(2+(2-p)|e|^2)+q\big)
=\w^2\w_q Q(|e|^2),
$$
where:
\begin{equation}
\label{eqQ}
Q(s)=Q_{p,q}(s)=2p+q+(2p+2q-p^2)s+qs^2.
\end{equation}
Thus the sign of $B$ is determined by the sign over $I_q$ of $Q(s)$.  If $q=0$ then $Q(s)=2p+p(2-p)s$, and $Q(I_q)>0$ precisely when $0<p\leqs2$.  If $q\neq0$ then $Q(s)$ has discriminant:
$$
E=p^2(p^2-4p-4q+4),
$$
and $E<0$ if and only if $p\neq0$ and $q>\k_1(p)$, where $\k_1(p)=\frac14(p-2)^2$.  Furthermore the critical point $s_0$ of $Q(s)$ satisfies $s_0<0$ if and only if $q>0$ and $q>\k_2(p)$, or $q<0$ and $q<\k_2(p)$, where $\k_2(p)=\frac12p(p-2)$.  It is helpful to note that:
\begin{align*}
0<\k_1(p)<\k_2(p)
\;\text{ and }\,
\k_1(p)<\lambda(p),
\quad
&\text{if $-8<p<-2$,} \\
0<\k_2(p)<\k_1(p),
\quad
&\text{if $-2<p<0$,} \\
\k_2(p)<0<\k_1(p),
\quad
&\text{if $0<p<2$,} \\
0<\k_1(p)<\k_2(p),
\quad
&\text{if $p>2$.}
\end{align*}
Suppose first that $q>0$.  Then $Q(s_0)$ is a global minimum.  If $q>\k_2(p)$ (ie. $s_0<0$) then $Q(\Bbb R^+)>0$ precisely when $2p+q>0$, whereas if $q\leqs\k_2(p)$ (ie. $s_0\geqs0$) then $|p|>2$ and hence $Q(\Bbb R^+)>0$ precisely when $q>\k_1(p)$ (ie. $E<0$).  
Now let $\Omega_i$ ($i=1,2,3$) be the following region of the half plane $2p+q>0$:
\begin{align*}
\Omega_1
&= \{(p,q): |p|>2 \text{ and } q>\k_1(p)\}, \\
\Omega_2
&= \{(p,q): -2\leqs p\leqs 0 \text{ and } 2p+q>0\}, \\
\Omega_3
&= \{(p,q): 0\leqs p\leqs 2 \text{ and } q>0\},
\end{align*}
and define:
$$
\Omega=\Omega_1\cup\Omega_2\cup\Omega_3.
$$
Thus $Q(I_q)>0$ if and only if $(p,q)\in\Omega$.  However $\Gamma_+\subset\Omega$, so no additional conditions are generated.  
Suppose now that $q<0$.  It suffices to consider $p>0$.  Now $Q(s_0)$ is a global maximum; and $E>0$, since $q<\k_1(p)$ and $p\neq0$.  Therefore, since $Q(0)=2p+q>0$, $Q(I_q)>0$ precisely when the maximum root $s_+\geqs-1/q$:
$$
\sqrt E\geqs p^2-2p-2q+2,
$$ 
which rearranges to:
$$
(p+q-1)^2\leqs0,
$$
and is therefore equivalent to $p+q=1$. 
In summary, additional conditions necessary for $\tilde K(\Pi)>0$ are $0<p\leqs2$ when $q=0$, and $p+q=1$ when $q<0$.

Finally we note that the necessary conditions listed above are in fact sufficient.  For, by \eqref{eq8} the sign of $\tilde K(X^v\wedge Y^v)$ is the same as that of $A ( \<X,e\>^2 + \<Y,e\>^2 ) + B$.  Now $0\leqs\<X,e\>^2 + \<Y,e\>^2\leqs|e|^2$, since $X$ and $Y$ are orthonormal.  But if $a,b\in\Bbb R$ with $a+b\geqs0$ and $b\geqs0$ then $au+b\geqs0$ for all $u\in[0,1]$.
\end{proof}

It is interesting to compare Theorem \ref{prop3a} with Corollary \ref{cor1}.
In \S3 we will extend Theorem \ref{prop3a} to a characterisation of $\tilde K\geqs0$ when $M$ has constant (non-negative) curvature (Theorem \ref{prop7}). 

We now turn to the Ricci curvature.  In the following result, it is understood that summation applies to all repeated indices.

\begin{proposition}\label{prop4}
Let $\rho$ denote the Ricci tensor of $(M,g)$.  Then the Ricci curvature $\tilde\rho$ of $h_{p,q}$ satisfies, at $e\in BM_q$:
\begin{align}
\tilde \rho(X^h,Y^h)
&= \rho(X,Y) - \tfrac34 \w^p \< R(X,e_i)e, R(Y,e_i)e \> \label{ric1} \\
&\qquad
+ \tfrac14 \w^p \< R(e,e_i)X, R(e,e_i)Y \> , \notag \\
\tilde \rho(X^h,Y^v)
&= \tfrac12 \w^p \< \d R(X)e,Y \>, \label{ric2} \\
\tilde \rho(X^v,Y^v)
&= \tfrac14 \w^{2p} \< R(X,e),R(Y,e) \> 
+ \a\<X,Y\> + \b\<X,e\>\<Y,e\>, \label{ric3}
\end{align}
for all vectors $X,Y \in T_{\pi(e)}M$, where  $\{e_i\}_{1\leqs i\leqs n}$ is an orthonormal basis of $T_{\pi(e)} M$, and $\d$ is the covariant coderivative of $(M,g)$.  The functions $\a$ and $\b$ on $BM_q$ depend on $A$ and $B$ of Proposition~\ref{prop2}:
$$
\a= |e|^2 \w_q A + (n-2+\w_q)B,
\qquad
\b= (n-1-\w_q)A + q\w_q B.
$$
In particular, $\a$ and $\b$ are independent of the curvature of $(M,g)$.
\end{proposition}

\begin{proof} 

Assume first that $e\neq0$, and let $\{e_1, ..., e_n\}$ be an orthonormal basis of $T_{\pi(e)} M$ with $e_1=e/|e|$.  Then $\{e_1^h, ..., e_n^h, f_1^v, ..., f_n^v\}$ is an orthonormal basis of $T_eTM$, where:
$$
f_1 = \sqrt{\w_q}\,\w^{-p/2}e_1 ,
\qquad 
f_j = \w^{-p/2}e_j ,
\quad\text{for $j\geqs 2$.}
$$
To prove formula \eqref{ric1}, we combine the definition of Ricci curvature,  parts (i) and (iv) of Proposition~\ref{prop2} and the definition of $h_{p,q}$ to obtain:
\begin{align*}
\tilde \rho(X^h,Y^h)
&= h_{p,q}( \tilde R(X^h, e_i^h)e_i^h, Y^h )
+ h_{p,q}( \tilde R(X^h, f_i^v)f_i^v, Y^h )  \\
&= \rho(X,Y) + \tfrac34 \w^p \< R(e, R(X,e_i)e)e_i, Y\>  
-\tfrac14 \w^p \< R(e,e_i) R(e,e_i)X, Y\>.
\end{align*}
To show \eqref{ric2}, we again use parts (i) and (iv) of Proposition~\ref{prop2} and the definition of $h_{p,q}$ to obtain:
\begin{align*}
\tilde \rho(X^h,Y^v)
&= h_{p,q}( \tilde R(X^h, e_i^h)e_i^h, Y^v)
+ h_{p,q}( \tilde R (X^h, f_i^v)f_i^v, Y^v)  \\
&=\tfrac12 \w^p \big[\big\<(\nabla_{e_i}R)(X, e_i)e, Y\big\> 
+ q\<Y,e\>\big\<(\nabla_{e_i}R)(X, e_i)e, e\big\>\big],
\end{align*}
and observe that $\<(\nabla_{e_i}R)(X,e_i)e, e\>=0$ to conclude.
For equation~\eqref{ric3}, we have:
\begin{align}
\label{eq*}
\tilde \rho(X^v, Y^v)
&= h_{p,q}( \tilde R(X^v, e_i^h)e_i^h, Y^v) 
+ h_{p,q}( \tilde R(X^v, f_i^v)f_i^v, Y^v).
\end{align}
First, using Proposition~\ref{prop2}\,(iv) and the definition of $h_{p,q}$, we obtain:
\begin{align*}
h_{p,q}( \tilde R(X^v, e_i^h)e_i^h, Y^v)
&= h_{p,q}( \tilde R(e_i^h,X^v)Y^v, e_i^h) \\
&= \tfrac14 \w^{2p} \< R(e,X)e_i, R(e,Y)e_i \>.
\end{align*}
We now use Proposition 2.3\,(vi) to expand the second term of \eqref{eq*}:
$$
\tilde R(X^v, f_i^v)f_i^v 
= A\,\<f_i,e\> [ r(X,f_i)e ]^v 
+ B [ r(X,f_i)f_i ]^v 
+ C \< r(X,f_i)f_i,e \> U,
$$
and then use $r(X,e_i)e_i=(n-1)X$, the relation $q|e|^2\w_q=1-\w_q$, and equation \eqref{eqABC} to rewrite this as:
\begin{align*}
\w^p \tilde R(X^v, f_i^v)f_i^v 
&= \big( |e|^2 \w_q A + (n-2+\w_q)B \big) X^v 
+ (n-2)C \<X,e\>U \\
&=\a X^v + (n-2)C\<X,e\>U.
\end{align*}
Therefore by the definition of $h_{p,q}$:
\begin{align*}
h_{p,q}( \tilde R(X^v, f_i^v)f_i^v, Y^v)
& =\big( (1-\w_q)A + q(n-2+\w_q)B +(n-2)\w_q^{-1} C \big)
\<X,e\>\<Y,e\> \\
&\qquad 
+ \a \<X,Y\> \\
& = \a\<X,Y\>+\b\<X,e\>\<Y,e\>,
\end{align*}
after applying equation \eqref{eqABC}.
This concludes the case $e\neq0$.
The formulas extend to $e=0$ by continuity.
\end{proof}

Recall that a Riemannian manifold is said to have {\sl harmonic curvature\/} if $\d R=0$.

\begin{corollary}\label{cor3}
Necessary conditions for $h_{p,q}$ to be an Einstein metric on $BM_q$ are that $(M,g)$ is Einstein and has harmonic curvature.
\end{corollary}

We also note from Proposition \ref{prop4}, setting $e=0$, that:
$$
\tilde\rho(X^h,Y^h)=\rho(X,Y),
\qquad
\tilde\rho(X^v,Y^v)=\a\<X,Y\>=(n-2)(2p+q)\<X,Y\>.
$$
Therefore necessary conditions for $\tilde\rho\geqs0$ are $\rho\geqs0$ and, when $n\geqs3$, $2p+q\geqs0$.  However, unlike the sectional curvature $\tilde K$, there are no values of $(p,q)$ for which $\tilde\rho$ extends to the tangent bundle of $SM_q$ (the problem being that $\a$ is always singular on $SM_q$).  

We conclude the section with a computation of the scalar curvature of $h_{p,q}$.

\begin{proposition}\label{prop5}
Let $s$ denote the scalar curvature of $(M,g)$.  Then for each $e\in BM_q$ the
scalar curvature $\tilde s$ of $h_{p,q}$ is:
$$
\tilde s
= s-\tfrac14 \w^p\sum_{i,j=1}^n |R(e_i,e_j)e|^2+(n-1)\w^{-p}(2\a-(n-2)B),
$$
where $\a$ and $B$ are as in Propositions~\ref{prop4} and \ref{prop2} respectively.
\end{proposition}

\begin{proof}
For $e\neq 0$, let $\{e_1^h, ..., e_n^h, f_1^v, ..., f_n^v\}$ be an orthonormal basis of $T_eTM$ as in the proof of Proposition~\ref{prop4}. By definition (summing over repeated indices):
$$
\tilde s
=  \tilde\rho (e_i^h, e_i^h)
+ \tilde\rho (f_i^v, f_i^v),
$$
and by Proposition \ref{prop4}:
\begin{align*}
\tilde\rho(e_i^h, e_i^h)
&= s - \tfrac34\w^p \sum_{i,j} |R(e_i,e_j)e|^2
+\tfrac14\w^p \sum_{i,j} |R(e,e_i)e_j|^2 \\
&= s - \tfrac14\w^p \sum_{i,j} |R(e_i,e_j)e|^2,
\end{align*} 
by the symmetries of the Riemann tensor, and:
$$
\tilde\rho(f_i^v, f_i^v)
= \tfrac14\w^p \sum_{i,j} |R(e,e_i)e_j|^2
+ \w^{-p}\big( (n-1+\w_q)\a + |e|^2\w_q\b\big).
$$
The result follows on noting that:
\begin{align*}
(n-1+\w_q)\a + |e|^2\w_q\b
&= (n-1)\big( 2|e|^2\w_q A + (n-2+2\w_q)B\big) \\
&=(n-1)(2\a - (n-2)B).
\end{align*}
For $e=0$, we use a continuity argument to finish the proof.
\end{proof}

\begin{remark}\label{rem3}
Setting $e=0$ in Proposition \ref{prop5} yields:
$$
\tilde s=s+(n-1)(2\a-(n-2)B)
=s+n(n-1)(2p+q).
$$
Therefore a necessary condition for $\tilde s>0$ is: 
$$
s>n(1-n)(2p+q),
$$ 
which does not preclude the possibility of $s<0$.  Note also that, as expected, $\tilde s$ does not extend smoothly across $SM_q$, for any values of $(p,q)$.
\end{remark}

\section{Generalised Cheeger-Gromoll metrics over space forms}

Unless otherwise stated, in this section $(M,g)$ is now a Riemannian manifold of dimension $n\geqs2$ and constant sectional curvature $c$.  The expressions for the Ricci, sectional and scalar curvatures of $h_{p,q}$ then simplify and their signs can be studied.

\begin{proposition}\label{prop6}
The Ricci, sectional and scalar curvatures of $h_{p,q}$ at $e\in BM_q$ are given by:
\begin{align}
\tilde \rho(X^h,Y^h) 
&= c(n-1)\<X,Y\> 
+ \tfrac12 c^2 \w^p \big[ (2-n) \<X,e\>\<Y,e\> 
- |e|^2 \< X,Y \> \big],
\label{eqRichh} \\
\tilde \rho(X^h,Y^v) &= 0,
\label{eqRichv} \\
\tilde \rho(X^v,Y^v) 
&= \big(\a + \tfrac12 c^2 |e|^2 \w^{2p} \big) \<X,Y\> 
+ \big( \b - \tfrac12 c^2 \w^{2p} \big) \<X,e\>\<Y,e\>,
\label{eqRicvv}
\end{align}
where $\a,\b$ are as defined in Proposition~\ref{prop4},
\begin{align}
\tilde K(X^h \wedge Y^h) 
&= c - \tfrac34 c^2 \w^p ( \<X,e\>^2 + \<Y,e\>^2 ),
\label{eqKhh}  \\
\tilde K(X^h \wedge Y^v) 
&= \tfrac14 c^2 \w^p \frac{\<X,e\>^2}{1+q\<Y,e\>^2 }, 
\label{eqKhv}  \\
\tilde K(X^v \wedge Y^v) 
&= \frac{\w^{-p}}{1 + q\<X,e\>^2 + q\<Y,e\>^2}
\big( A (\<X,e\>^2 + \<Y,e\>^2 ) +B \big),
\label{eqKvv}
\end{align}
for all orthonormal vectors $X$ and $Y$ of $T_{\pi(e)} M$, 
\begin{equation}
\tilde s = (n-1)\big( nc - \tfrac12 c^2 \w^p |e|^2 
+ \w^{-p} (2\a - (n-2)B)\big). 
\label{eqscal}
\end{equation}
\end{proposition}

\begin{proof}
In this situation $R=cr$, where $r$ is as defined in \eqref{littler}.
We first compute the sectional curvatures.  We have, for $X$ and $Y$ orthonormal:
$$ 
|R(X,Y)e|^2 = c^2 ( \<X,e\>^2 + \<Y,e\>^2),
$$
so formula~\eqref{eq6} of Proposition~\ref{prop3} yields \eqref{eqKhh}.
Similarly, $R(e,Y)X = -c\<X,e\>Y$ and formula~\eqref{eq7} of Proposition~\ref{prop3} implies \eqref{eqKhv}.  Formula (8) for $\tilde K(X^v \wedge Y^v)$ is unchanged, being independent of the curvature of $(M,g)$.

For the scalar and Ricci curvatures, let $\{e_1,\dots,e_n\}$ be an orthonormal basis of $T_{\pi(e)}M$ with $e_1=e/|e|$ (assuming $e\neq 0$).  Then:
$$
\sum_{i,j=1}^n |R(e_i,e_j)e|^2 
= c^2 |e|^2 \sum_{i\neq j}^n \big( \<e_i,e_1\>^2 + \<e_j,e_1\>^2 \big)
= 2(n-1) c^2 |e|^2.
$$
Furthermore $s=n(n-1)c$, so Proposition~\ref{prop5} for $\tilde s$ reduces to formula~\eqref{eqscal}.  The case $e=0$ follows by continuity.  For the Ricci curvature, we note first that (summing over repeated indices):
\begin{align*}
\<R(X,e_i)e,R(Y,e_i)e\> 
&= c^2 \big( |e|^2\<X,Y\> + (n-2)\<X,e\>\<Y,e\> \big) \\
&= \<R(e,e_i)X,R(e,e_i)Y\>. 
\end{align*}
Furthermore $\rho(X,Y) = c(n-1)\<X,Y\> $.  Plugging these into formula~\eqref{ric1} yields \eqref{eqRicvv}.  Equation (14) follows from (10) and the fact that space forms have harmonic curvature.  Finally:
$$
\<R(X,e)e_i,R(Y,e)e_i\> 
= 2c^2 \big( |e|^2 \<X,Y\> - \<X,e\>\<Y,e\> \big),
$$
from which formula (11) reduces to (15).  (Note that $\a$ and $\b$ are independent of the curvature of $M$.)  The case $e=0$ again follows by continuity.
\end{proof}

Proposition \ref{prop6} allows us to investigate conditions under which the sectional or scalar curvature of $h_{p,q}$ is non-negative.  Note from equation \eqref{eqKhv} that the sectional curvature is never strictly positive.  In fact we reach a characterisation of $\tilde K \geqs 0$, and sufficient conditions for $\tilde s>0$.  Recall that a necessary condition for $\tilde K\geqs 0$ is $K\geqs0$, and hence $c\geqs0$.  
Recalling the functions $\mu(p)$ and $\lambda(p)$, defined in equations \eqref{eqmu} and \eqref{eqlambda} respectively, we now modify the regions $\Gamma$ and $\Gamma'$ introduced in Definition \ref{def1} to the following families of regions $\Delta_c$ and $\Delta_c'$, for $c\geqs0$.

\begin{definition}\label{def2}
Let $\Delta_+^i$ ($i=1,2,3$) be the following region of the $(p,q)$-plane:
\begin{align*}
\Delta_+^1
&= \{(p,q): -8<p\leqs -2 \text{ and } q\geqs\lambda(p)\}
\supset\Gamma_+^1, \\
\Delta_+^2
&= \{(p,q): -2\leqs p\leqs 0 \text{ and } 2p+q\geqs0\}
\supset\Gamma_+^2, \\ 
\Delta_+^3
&= \{(p,q): 0\leqs p\leqs 1 \text{ and } q>0\}
=\Gamma_+^3,
\end{align*}
and set:
$$
\Delta_+=\Delta^1_+\cup\Delta^2_+\cup\Delta^3_+.
$$
Define further regions:
\begin{align*}
\Delta_-&=\{(p,q):p+q=1 \text{ and } q<0\}=\Gamma_-, \\
\Delta'_-&=\{(p,q):p+q\geqs1 \text{ and } q<0\}=\Gamma'_-, \\
\Delta_Z&=\{(p,0):0\leqs p\leqs2\}\supset\Gamma_Z, \\
\Delta'_Z&=\{(p,0):p\geqs0\}\supset\Gamma'_Z,
\end{align*}
and set:
$$
\Delta_0=\Delta_-\cup\Delta_Z\cup\Delta_+,
\qquad
\Delta'_0=\Delta'_-\cup\Delta'_Z\cup\Delta_+.
$$  
Notice that $\Delta_0$ (respectively, $\Delta_0'$) lies in the closure of $\Gamma$ (respectively, $\Gamma'$).  Now for $c>0$ define:
\begin{align*}
\Delta_{c,-}&=\{(p,q):p+q=1,\ q<0 \text{ and } \mu(p)\geqs3c/4\}
\subset\Delta_-, \\
\Delta'_{c,-}&=\{(p,q):p+q\geqs1,\ q<0 \text{ and } \mu(p)\geqs3c/4\}
\subset\Delta'_-, \\
\Delta_{c,Z}&=\{(p,0):1\leqs p\leqs2 \text{ and } \mu(p)\geqs3c/4\}
\subset\Delta_Z, \\
\Delta_{c,Z}'&=\{(p,0):p\geqs1 \text{ and } \mu(p)\geqs3c/4\}
\subset\Delta_Z', \\
\Delta_{c,+}&=
\left\lbrace
\begin{array}{l}
\{(1,q):q>0\},
\quad\text{if $c\leqs 4/3$} \\
\emptyset,
\quad\text{if $c>4/3$}
\end{array}
\right\rbrace
\subset\Delta_+,
\end{align*}
and set: 
$$
\Delta_c=\Delta_{c,-}\cup\Delta_{c,Z}\cup\Delta_{c,+},
\qquad
\Delta_c'=\Delta'_{c,-}\cup\Delta_{c,Z}'\cup\Delta_{c,+}.
$$
Notice that $\Delta_{c,Z}=\emptyset$ if $c>4\mu(2)/3=16/3$.
\end{definition}

For all $c\geqs0$ we have $\Delta_c\subset\Delta_c'$, and $\Delta_c'$ lies in the closed half plane $2p+q\geqs0$.  In fact, if $c>0$ then $\Delta_c'$ lies in the region $p+q\geqs1$ with $p\geqs1$, and is independent of $c$ if $c\leqs4/3$, whereas if $c>16/3$ then $\Delta_c'$ lies in the region $p+q\geqs1$ with $q\leqs0$ and $p>2$.  Furthermore $\Delta_{c_1}\subset\Delta_{c_2}\subset\Gamma\subset\Delta_0$ for all $c_1>c_2>0$, and similar relations hold for the $\Delta_c'$.  Note that $\Delta_-$ (respectively, $\Delta'_-$) is the limit of $\Delta_{c,-}$ (respectively, $\Delta'_{c,-}$) as $c\to0$, but the corresponding limits do not hold for the other components of $\Delta_c$ and $\Delta'_c$.  In particular, $\Delta_c$ and $\Delta'_c$ do not converge to $\Delta_0$ and $\Delta'_0$ respectively, as $c\to0$.

\begin{theorem}\label{prop7}
Let $M$ be an $n$-dimensional Riemannian manifold of constant sectional curvature $c\geqs0$.  If $n\geqs3$ {\rm(}respectively, $n=2${\rm)} then $h_{p,q}$ has non-negative sectional curvature precisely when $(p,q)\in\Delta_c$ {\rm(}respectively, $(p,q)\in\Delta_c'${\rm)}.
\end{theorem}

\begin{proof}
Necessary and sufficient conditions for $\tilde K(\Pi)\geqs0$ for all vertical $2$-planes $\Pi$ may be deduced from Theorem \ref{prop3a}, yielding the regions $\Delta_0$ and $\Delta_0'$.  When $c=0$ it follows from Proposition~\ref{prop6} that these are necessary and sufficient for $\tilde K\geqs0$.   Suppose $c>0$.  It follows from equation \eqref{eqKhv} that $\tilde K(\Pi)\geqs0$ for all vertizontal $\Pi$.  For all dimensions,  it follows from \eqref{eqKhh} that $\tilde K(\Pi) \geqs 0$ for all horizontal $2$-planes $\Pi$ in $T_eTM$ if and only if:
\begin{equation}
\w^p |e|^2 \leqs \frac{4}{3c}.
\label{ineq}
\end{equation}
Now $\w^p|e|^2=f(|e|^2)$ where $f(t)=t/(1+t)^p$.  The function $f(t)$ is bounded on $\Bbb R^+$ if and only if $p\geqs1$, in which case it has supremum $1/\mu(p)$, attained at $t=1/(p-1)$ if $p>1$.  If $q\geqs0$ then it follows that $\tilde K(\Pi)\geqs 0$ for all horizontal $\Pi$ precisely when $p\geqs 1$ and $\mu(p)\geqs 3c/4$.  Imposing these conditions on $\Delta_+$, $\Delta_Z$  and $\Delta_Z'$ yields the regions $\Delta_{c,+}$, $\Delta_{c,Z}$ and $\Delta'_{c,Z}$, respectively.  If $q<0$ then we only require \eqref{ineq} to hold for $|e|^2<-1/q$.  Now if $p>1$ then $\w^p|e|^2$ has supremum $1/\mu(p)$ over $BM_q$ precisely when $p+q\geqs1$.  
Since $p>1$ and $p+q\geqs1$ are both necessary conditions (from $\Delta_-$ and $\Delta'_-$), it follows that \eqref{ineq} holds if and only if $\mu(p)\geqs 3c/4$.
\end{proof}

It follows from Theorem \ref{prop7} that if $0\leqs c\leqs 4/3$ then the conditions for $\tilde K\geqs0$ are independent of $c$,
whereas if $n\geqs3$ (repectively, $n=2$) and $c>16/3$ then all $(p,q)$ for which $h_{p,q}$ has non-negative sectional curvature lie in the open (respectively, closed) lower half plane $q<0$ (respectively, $q\leqs0$).

We now consider the scalar curvature $\tilde s=\tilde s_{p,q}$ of $h_{p,q}$.  Since $s=n(n-1)c$ it follows from Remark \ref{rem3} that a necessary condition for $\tilde s>0$ is $2p+q>-c$.  Notice that this no longer precludes $c<0$.  We will not attempt to find all metrics with $\tilde s>0$, but identify subregions of $\Gamma$ and $\Gamma'$ where this is the case.  It will be convenient to introduce the following function:
\begin{equation}
\label{eqnu}
\nu(p)=
\dfrac{2(1-p)}{2+p},
\quad
p\neq-2,
\end{equation}
which parametrises the hyperbola in the $(p,q)$-plane with equation $pq+2p+2q=2$.  Notice that $1-p<\lambda(p)<\nu(p)<0$ for all $p>1$.  We now introduce the following ``multipliers'' $m_i=m_i(n,p,q)\geqs1$ ($i=1,\dots,5$).

\begin{definition}\label{def3}
If $p>0$ define $m_1\geqs1$ by:
$$
(m_1)^2=1+\frac{2(n-2)}{n^2p}.
$$
Note that $m_1=1$ if and only if $n=2$.  
If $p\geqs1$ define $m_2>1$ and $m_3\geqs1$ by:
$$
(m_2)^2=1+\frac{4p}{n\mu(p)},
\qquad
(m_3)^2=1+\frac{2(p^2-1)}{np\mu(p)}.
$$
Note that $m_3>1$ if $p>1$.
If $p>1$ and $\lambda(p)<q<0$ define $m_4>1$ by:
$$
(m_4)^2=1+\frac{D}{4(p-1)(q-1)\mu(p)},
$$
where:
$$
D=pq(pq+8p+8q-8).
$$
Note that $D<0$ for all $(p,q)$ in the region $\lambda(p)<q<0$.  
Finally, if $p>1$ and $p+q\geqs1$ define $m_5\geqs1$ by:
$$
(m_5)^2=1+\frac{p+q-1}{\mu(p)}.
$$
Note that $m_5=1$ if and only if $p+q=1$, and $m_5=m_4$ when $q=\nu(p)$.
\end{definition}

\begin{theorem}\label{prop8}
Let $(M,g)$ be an $n$-dimensional Riemannian manifold of constant sectional curvature $c$, and let $\tilde s$ be the scalar curvature of $h_{p,q}$ on $BM_q$. 

Let $n=2$.
For $c=0$, $\tilde s>0$ precisely when $(p,q)\in\Gamma'$.   
For $c\neq0$, $\tilde s>0$ if one of the following holds: \\
{\rm (a)}\quad 
$q>0$, $p=1,$ and $|c-2|<2$, \\
{\rm (b)}\quad
$q=0$, $1\leqs p<2$, and $|c-2\mu(p)|<2\mu(p)$, \\
{\rm (c)}\quad 
$q=0$, $p\geqs2$, and $|c-2\mu(p)|<2m_2\mu(p)$, \\
{\rm (d)}\quad 
$q<0$, $p>1$, $q\geqs\nu(p)$, and $|c-2\mu(p)|<2m_4\mu(p)$, \\
{\rm (e)}\quad 
$q<0$, $p+q\geqs1$, $q\leqs\nu(p)$, and $|c-2\mu(p)|<2m_5\mu(p)$. 

\smallskip
Let $n\geqs 3$.
For $c=0$, $\tilde s>0$ if $(p,q)\in\Gamma$.
For $c\neq 0$, $\tilde s>0$ if one of the following conditions holds: \\
{\rm(a)}\quad 
$q>0$, $p=1$, and $|c-n|<m_1n$, \\
{\rm(b)}\quad
$q=0$, $1\leqs p<2$, and $|c-n\mu(p)|<m_1n\mu(p)$, \\
{\rm (c)}\quad 
$q=0$, $p=2$, and $|c-4n|<4m_2n$, \\
{\rm (d)}\quad 
$q<0$, $p+q=1$, and $|c-n\mu(p)|<m_3n\mu(p)$.

\comment
If $n\geqs 3$ then: \\
{\rm(1)}\quad
For $c=0$, $\tilde s>0$ if $(p,q)\in\Gamma$, and also in the following cases:

{\rm i)}\quad 
if $q=0$ then $\tilde s \geqs 0$ precisely when $0< p\leqs 2$,

{\rm ii)}\quad 
if $q>0$ and $p=0,1$ or $2$ then $\tilde s \geqs 0$,

{\rm iii)}\quad 
if $q>0$, $0<p<2$ and $n\geqs4$ then $\tilde s \geqs 0$. \\
{\rm(2)}\quad 
For $c \neq 0$:

{\rm i)}\quad 
if $q=0$ and $p=1$, $\tilde s>0$ precisely when $|c-n| < \sqrt{n^2 + 2n-4}$,
    
{\rm ii)}\quad 
if $q=0$ and $p=2$, $\tilde s>0$ precisely when $|\frac{c}{4}-n| < \sqrt{n(n+2)}$,
    
{\rm iii)}\quad 
if $q=0$, $1<p<2$ and $0<c<2n\mu(p)$ then $\tilde s> 0$.
\endcomment

\end{theorem}

\begin{proof}
From Proposition \ref{prop6}, after some calculation:
\begin{equation}
\label{eqcscal}
\frac{\tilde s}{n-1}=nc-\tfrac12 c^2f(t)+\v(t),
\end{equation}
where $t=|e|^2$.  Here $f(t)=t/(1+t)^p$, the function introduced in the proof of Theorem \ref{prop7}, and:
\begin{equation}
\label{eqphi}
\v(t)=(1+t)^{p-2}(1+qt)^{-2}C(t),
\end{equation}
with $C(t)$ the following generically cubic polynomial:
\begin{equation}
\label{eqC}
C(t)=C_{p,q}(t)=2P(t)+(n-2)(1+qt)Q(t),
\end{equation}
involving the polynomials $P(t)$ and $Q(t)$ defined in \eqref{eqP} and \eqref{eqQ}, respectively.  If $c=0$ then $\tilde s>0$ precisely when $\v(I_q)>0$, where $I_q$ is the interval defined in \eqref{eqIq}.  When $n=2$,  $\v(I_q)>0$ if and only if $P(I_q)>0$, and from the proof of Theorem \ref{prop3a} this is the case precisely when $(p,q)\in\Gamma'$.  When $n\geqs3$, a sufficient condition for $\v(I_q)>0$ is $P(I_q)>0$ and $Q(I_q)>0$, and from the proof of Theorem \ref{prop3a} this is the case precisely when $(p,q)\in\Gamma$.
If $c\neq0$ we restrict attention to those $(p,q)\in\Gamma$ (if $n\geqs3$), or $(p,q)\in\Gamma'$ (if $n=2$), with $p\geqs1$.  This ensures that $f(t)$ is bounded on $\Bbb R^+$, with supremum $1/\mu(p)$, and $\v(t)$ is bounded below on $I_q$, with infimum $\check\v\geqs0$ (Theorem \ref{prop3a}).  Then \eqref{eqcscal} yields the following sufficient condition for $\tilde s>0$:
$$
\frac{c^2}{2\mu(p)}<nc+\check\v,
$$
which is equivalent to the following quadratic inequality for $c$:
\begin{equation}
\label{cineq}
c^2-2n\mu(p)c-2\mu(p)\check\v<0.
\end{equation}

Suppose first that $n=2$.  
Referring to Theorem \ref{prop3a}, the possibilities for $(p,q)\in\Gamma'$ with $p\geqs1$ are as follows. \\
(a)\quad
$q>0$ and $p=1$.
Then $\check\v=0$, and \eqref{cineq} yields $0<c<4$. \\
(b)\quad
$q=0$ and $1\leqs p<2$.
Then $\check\v=0$, and \eqref{cineq} yields $0<c<4\mu(p)$.  \\
(c)\quad
$q=0$ and $p\geqs2$.
Then $\check\v=4p$, and \eqref{cineq} yields:
$$
|c-2\mu(p)|<2\mu(p)\sqrt{1+2p/\mu(p)}=2m_2\mu(p).
$$
When $q<0$ we estimate $\check\v$ by first noting the following lower bounds over $I_q$:
\begin{align}
(1+qt)^{-2}&\geqs1,
\notag  \\
(1+t)^{p-2}&\geqs
\left\{
\begin{array}{l}
\Big(\dfrac{q-1}{q}\Big)^{p-2},\quad\text{if $1<p<2$} \\
1,\quad\text{if $p\geqs2$}
\end{array}
\right\}
\geqs\dfrac{q}{q-1}. 
\label{lbound}
\end{align}
Referring to the proof of Theorem \ref{prop3a}, $P(t)$ has a global minimum at $t_0>0$, and we note that $t_0\geqs-1/q$ if and only if $pq+2p+2q-2\leqs0$.  It follows that if $q\geqs\nu(p)$ then:
$$
\inf_{t\in I_q} P(t)=P(t_0)=\frac{D}{4(p-1)q},
$$
whereas if $q\leqs\nu(p)$ then:
$$
\inf_{t\in I_q}P(t)=P(-1/q)=\frac{(p+q-1)(q-1)}{q}.
$$
These estimates yield the following two further cases. \\
(d)\quad
$q<0$, $p>1$ and $q\geqs\nu(p)$.
Then $\check\v\geqs D/2(p-1)(q-1)$, and \eqref{cineq} is therefore satisfied if:
$$
|c-2\mu(p)|<2m_4\mu(p).
$$
(e)\quad
$q<0$, $p+q\geqs1$ and $q\leqs\nu(p)$.  Then $\check\v\geqs2(p+q-1)$, and \eqref{cineq} is satisfied if:
$$
|c-2\mu(p)|<2m_5\mu(p).
$$
This completes our analysis of surfaces.

Now suppose that $n\geqs 3$.  Referring to Theorem \ref{prop3a}, the possibilities for $(p,q)\in\Gamma$ with $p\geqs1$ are as follows. \\
(a)\quad
$q>0$ and $p=1$.
Then $\v$ is a rational function, which is smooth and decreasing on $\Bbb R^+$, hence:
$$
\check\v=\lim_{t\to\infty}\v(t)=n-2.
$$
Therefore \eqref{cineq} is equivalent to: 
$$
|c-n|<n\sqrt{1+2(n-2)/n^2}=m_1n.
$$
(b)\quad
$q=0$ and $1\leqs p<2$.  Then:
$$
\v(t)=p(1+t)^{p-2}(2n+(n-2)(2-p)t).
$$
If $p=1$ then $\check\v=n-2$, and if $1<p<2$ then $\v$ has a global minimum at: 
$$
\tau=(n+2)/(n-2)(p-1),
$$ 
hence:
$$
\check\v
=\v(\tau)
=(n-2)^{2-p}(p-1)^{1-p}((n-2)p+4)^{p-1}
>(n-2)\mu(p)/p.
$$
It follows that $\check\v\geqs(n-2)\mu(p)/p$ for all $1\leqs p<2$.  Therefore \eqref{cineq} is satisfied if:
$$
|c-n\mu(p)|<n\mu(p)\sqrt{1+2(n-2)/n^2p}=m_1n\mu(p).
$$
(c)\quad
$q=0$ and $p=2$.
Then $\check\v=4n$ and $\mu(p)=4$, so \eqref{cineq} is equivalent to:
$$
|c-4n|<4n\sqrt{1+2/n}=4nm_2.
$$
(d)\quad
$q<0$ and $p+q=1$.  Then:
$$
P(t)=(1+qt)(2-q+qt),
\quad
Q(t)=(1+qt)(2-q+t),
$$
and it follows that:
$$
\v(t)=(1+t)^{p-2}\Big((n-2)(2-q+t)+\frac{2(2-q+qt)}{1+qt}\Big)
=(1+t)^{p-2}\psi(t),
\quad\text{say.}
$$
Now $\psi(t)$ is increasing on $I_q$, and therefore bounded below by $\psi(0)=n(2-q)$.  Also $(1+t)^{p-2}$ is bounded below on $I_q$ by $(p-1)/p$, from \eqref{lbound}.  Therefore $\check\v\geqs n(p^2-1)/p$, so \eqref{cineq} is satisfied if: 
$$
|c-n\mu(p)|<n\mu(p)\sqrt{1+2(p^2-1)/np\mu(p)}=m_3n\mu(p).
$$
\end{proof}

As a correction to Sekizawa's computations of the sectional and scalar curvatures of the Cheeger-Gromoll metric $h_{1,1}$, Gudmundsson and Kappos proved the following result.

\begin{theorem}\label{thmgudkap}\cite{Gud-Kap2}
Let $(M,g)$ be an $n$-dimensional Riemannian manifold of constant sectional curvature $c$, then there exist real numbers $c_n$ and $C_n$ such that $(TM,h_{1,1})$ has positive scalar curvature if and only if $c\in (c_{n},C_{n})$.  If $n=2$ then $c_2=0$ and $C_2 \geqs 40$, and if $n\geqs 3$ then $c_n<0$ and $C_n>60$.
\end{theorem}

Our next result shows that the curvature restrictions of Theorem \ref{thmgudkap} can be lifted by varying the parameters $(p,q)$.   

\begin{theorem}\label{thm1}
Let $(M,g)$ be an $n$-dimensional Riemannian manifold {\rm ($n\geqs 2$)} of constant sectional curvature $c$. Then there exist parameters $p$ and $q$ such that $h_{p,q}$ has positive scalar curvature.
\end{theorem}

\begin{proof}
If $n=2$ then it follows from Theorem \ref{prop8} that if $c=0$ (respectively, $c>0$) then $h_{p,0}$ has positive scalar curvature for all $p>0$ (respectively, all $p\geqs2$ with $\mu(p)>c$), whereas if $c<0$ then $\tilde s_{p,0}>$ provided $p\geqs2$ and:
$$
2(m_2-1)\mu(p)>-c.
$$
Note that:
$$
\frac{p}{\mu(p)}=\Big(\frac{p-1}{p}\Big)^{p-1}
$$
is decreasing, with limit $1/e$ as $p\to\infty$, so $m_2$ is bounded below by $\sqrt{1+2/e}$.  Therefore $\tilde s_{p,0}>0$ for all $p\geqs2$ satisfying:
$$
\mu(p)>\frac{-c}{\sqrt{e+2}-\sqrt e}.
$$

Suppose now that $n\geqs 3$, and consider $h_{p,q}$ with $p+q=1$ and $q<0$.  By Theorem \ref{prop8}, if $c=0$ (respectively, $c>0$) then $h_{p,q}$ has positive scalar curvature for all such $(p,q)$ (respectively, all such $(p,q)$ with $\mu(p)>c$), whereas if $c<0$ then $\tilde s_{p,q}>0$ provided:
$$
n(m_3-1)\mu(p)>-c.
$$
Now $m_3$ is increasing (in $p$), and hence bounded below by $\sqrt{1+3/4n}$ when $p\geqs2$.  Therefore $\tilde s_{p,q}>0$ if $p\geqs2$ and (for example):
$$
\mu(p)>\frac{-c}{\sqrt{n+3/4}-\sqrt n}.
$$
\end{proof}

One difficulty with Theorem \ref{thm1} is that when $n\geqs3$ the metrics $h_{p,q}$ with $\tilde s_{p,q}>0$ all have $q<0$, and therefore only endow the Riemannian ball subbundle of $TM$ with a metric of positive scalar curvature.
If $c=0$ then by Theorem \ref{prop8} any metric $h_{p,q}$ with $0\leqs p\leqs1$ and $q>0$ has positive scalar curvature, and our next result generalises this to all values of $c$.  Scrutiny of the proof will show, however, that when $n\geqs3$ and $c\neq0$ the parameters $(p,q)$ in general lie in the interior of the complement of $\Gamma$, which by Theorem \ref{prop3a} implies that the sectional curvatures of the tangent spaces are no longer entirely positive (or non-negative).

\begin{theorem}\label{thm3}
Let $(M,g)$ be an $n$-dimensional Riemannian manifold {\rm ($n\geqs 2$)} of constant sectional curvature $c$. Then there exist parameters $p$ and $q\geqs0$ such that $h_{p,q}$ has positive scalar curvature.
\end{theorem}

\begin{proof}
It suffices to consider $n\geqs3$ and $c\neq0$.
We expand the polynomial $C(t)$ defined in equation \eqref{eqC}:
$$
C(t)=(n-2)q^2t^3+a_{p,n}(q)t^2+b_{p,n}(q)t+n(2p+q),
$$
where:
\begin{align*}
a_{p,n}(q)
&=2(n-2)q^2+(n+2(n-3)p-(n-2)p^2)q, \\
b_{p,n}(q)
&=2(np-1)q^2+2(n+(n-1)p)q+(n-2)p(2-p).
\end{align*}
For fixed $p\geqs2$, the roots of $a_{p,n}(q)$ (respectively, $b_{p,n}(q)$) are real, and the larger root $A_{p,n}$ (respectively, $B_{p,n}$) is non-negative.  Therefore, if $q\geqs0$ then $C(t)$ will be positive on $\Bbb R^+$ if in addition $q>\max(A_{p,n},B_{p,n})$.  Thus $\check\v\geqs0$, and by \eqref{cineq} a sufficient condition for $\tilde s>0$ is:
$$
c^2-2n\mu(p)c<0.
$$
Therefore if $c>0$ then $\tilde s>0$ by choosing $\mu(p)>c/2n$. 

When $c<0$, we assume that $p$ is a positive integer.  It then follows from \eqref{eqcscal} that the sign of $\tilde s$ is controlled by the polynomial:
$$ 
G(t)
=nc(1+t)^{p}(1+qt)^2
-\tfrac12 c^2t(1+ qt)^2 
+(1+t)^{2p-2} C(t),
$$ 
for $t\geqs 0$.
We can make the coefficients of degree greater or equal to $p+3$ positive by
requiring that, for a chosen $p$, $q$ be greater than $q_{p+3} =\max(A_{p,n},B_{p,n})$.
The constant in $G(t)$ is  positive if, for $p$ fixed, $q$ is greater than $q_0 = -c - 2p$, and the coefficient of $t$ in $G(t)$ is a polynomial in $q$, of degree two and positive leading coefficient, so that it is positive if $q$ is greater than a certain real number $q_1$.  Similarly the coefficients of $t^2$ and $t^3$ are polynomials of degree two in $q$, with positive leading coefficients if $p\geqs p_2$ and $p\geqs p_3$, respectively, where $p_2,p_3>1$; so for $q$ greater than some real numbers $q_2$ and $q_3$, the coefficients of $t^2$ and $t^3$ in $G(t)$ are positive.  So if $p$ is greater than $\max(p_2,p_3)$ and $q$ is greater than $\max(q_0,q_1,q_2,q_3,q_{p+3})$, the constant term and the coefficients of $t$, $t^2$, $t^3$ and $t^k$ ($k\geqs p+3$) in $G(t)$, will be positive.

The coefficients of $t^k$ for $4\leqs k\leqs p+2$ are given by polynomials of degree two in $q$, with leading coefficients positive if:
\begin{equation}
\label{eq**}
p\geqs\frac{k-2}{2}-\frac{(k-1)nc}{(n-2)2^{k-1}} ,
\end{equation}
for all $ k=4,\dots,p+2$.
Studying the function
$$
\rho(x)=\frac{x-2}{2}-\frac{nc(x-1)}{(n-2)2^{x-1}}
$$
over the interval $[4,p+2]$, shows that:
$$ 
\max_{x\in[4,p+2]}{\rho(x)}=\max(\rho(4),\rho(p+2)).
$$
Using the upper bounds $n/n-2\leqs3$ and $(p+1)/p\leqs2$ for all $n\geqs3$ and all $p\geqs1$, it follows that \eqref{eq**} is satisfied if:
$$
p \geqs p_4=\max(1-9c/8,1+\ln(-3c)/\ln 2),
$$ 
and the leading terms of the coefficients of $t^k$, for $4\leqs k\leqs p+2$, will be positive.

To conclude the proof, choose $p\geqs\max(p_2,p_3,p_4)$, ensuring that all the leading terms of the coefficients of $G(t)$ (seen as second order polynomials in $q$) are positive, and then, for this fixed value of $p$, choose $q$ greater than:
$$
\max(q_0,q_1,q_2,q_3,q_4,\dots q_{p+2},q_{p+3}),
$$
where, for $4\leqs k\leqs p+2$, $q_k$ is the (larger) positive root of the coefficient of $t^k$ in $G(t)$. For such $p$ and $q$, $G(t)$ has all its coefficients positive, hence $\tilde s$ is positive.
\end{proof}

\begin{remark}
In the context of $g$-natural metrics, results close to Theorems~\ref{thm1} and \ref{thm3} can be found in~\cite{ab1} and~\cite{ab2}. More precisely, \cite[Theorem 1.8]{ab1} states that on the tangent bundle of a Riemannian manifold of negative scalar curvature, one can construct a $g$-natural metric of negative scalar curvature, while \cite[Theorem A.2]{ab2} shows that, on a manifold of constant negative  sectional curvature, one can find functions endowing the tangent space with a $g$-natural metric of  constant positive  scalar curvature.  However, it should be noted that even under the hypothesis of \cite[Theorem A.2]{ab2}, Theorem~\ref{thm3} identifies the positive scalar curvature metric on $TM$ as a simpler metric, belonging to a tighter class, which is a 2-parameter variation of the original Cheeger-Gromoll metric and for which the map $\pi\colon TM \to M$ is a Riemannian submersion with totally geodesic fibres.
\end{remark}

Theorem~\ref{thm1} can be extended to more general manifolds under some curvature conditions.

\begin{theorem}\label{prop9}
Let $(M,g)$ be an $n$-dimensional Riemannian manifold {\rm ($n\geqs 2$)}. If there exist constants $a$ and $b$ such that:
$$ 
s\geqs a 
\quad\text{and}\quad 
\sum_{i,j=1}^{n} |R(e_i,e_j)e|^2\leqs b|e|^2,
$$ 
for all $e\in TM$, then there exist parameters $(p,q)$ such that $(TM,h_{p,q})$ has positive scalar curvature.
\end{theorem}

\begin{proof}
By Proposition \ref{prop5}, for any $(p,q)$:
$$
\tilde s
\geqs a-\tfrac14 b\/\w^p|e|^2+(n-1)\w^{-p}(2\a-(n-2)B).
$$
It follows from Proposition \ref{prop6} that if $c$ is chosen to be less than: 
$$ 
\min\big(a/n(n-1),-\sqrt{b/2(n-1)}\,\big),
$$
then $\tilde s\geqs\tilde s(c)$, the scalar curvature of the generalised Cheeger-Gromoll metric on $BM_q$ if $(M,g)$ were a space form of curvature $c$.  By Theorem \ref{thm3}, it is possible to choose $(p,q)$ with $q\geqs0$  such that $\tilde s(c)>0$.
\end{proof}

The curvature conditions of Theorem \ref{prop9} are always satisfied on compact manifolds.

\begin{theorem}\label{thm2}
Let $(M,g)$ be a compact Riemannian manifold. Then there exist parameters $p$ and $q$ such that $(TM,h_{p,q})$ has positive scalar curvature.
\end{theorem}

\begin{proof}
Since $M$ is compact and the scalar curvature is a continuous function, the existence of the constant $a$ of Theorem \ref{prop9} is automatic. 

To establish the existence of a constant $b$, we will start from a finite open covering $\mathcal U$ of $M$ such that on (the closure of) each $U\in\mathcal U$ there exists a local orthonormal frame $\{e_1,\dots,e_n\}$.  To obtain $b$, for each $U$ and each pair $(e_i,e_j)$, we need only establish the existence of a constant $M_{ij}$ such that:
$$ 
|R(e_i,e_j)e|^2 \leqs M_{ij}|e|^2,
$$
for all $e\in T_xM$ and all $x\in U$.
If $e=\sum_{k=1}^{n} a_k e_k$ then: 
$$
|R(e_i,e_j)e|^2 = |\sum_k a_k R(e_i,e_j)e_k|^2 \leqs \sum_k |a_k|^2  |R(e_i,e_j)e_k|^2.
$$
For each $k$, the function $x \mapsto |R_{x}(e_i,e_j)e_k|_{x}^2$, defined on the relatively compact set $U$, is continuous and therefore has a supremum $M_{ijk}$.  Taking $M_{ij} = \max_k M_{ijk}$ yields:
$$
|R(e_i,e_j)e|^2 \leqs M_{ij} \sum_k |a_k|^2 = M_{ij} |e|^2 .
$$
\end{proof}

\end{document}